\documentclass[12pt]{article}
\usepackage{amssymb}
\newtheorem{proposition}{Proposition}[section]  
\newcommand{\bprop}{\medskip\begin{proposition}  \it}
\newcommand{\eprop}{\end{proposition} \hfill  \\}

\newcommand{\bed}{\begin{displaymath}}
\newcommand{\eed}{\end{displaymath}}
\newcommand{\ba}{\begin{array}}
\newcommand{\ea}{\end{array}}
\newcommand{\beq}{\begin{equation}}
\newcommand{\eeq}{\end{equation}}
\newcommand{\be}{\begin{equation}}
\newcommand{\ee}{\end{equation}}
\newcommand{\bea}{\begin{eqnarray}}
\newcommand{\eea}{\end{eqnarray}}

\def\lra{\longrightarrow}
\def\ot{{\otimes}}
\def\Z{{\mathbb Z}}
\def\N{{\mathbb N}}
\def\C{{\mathbb C}}

\def\M{\mathrm{Mat}}

\def\st{\stackrel}

\begin{document}
\baselineskip17.5pt
\setcounter{page}{0}
\thispagestyle{empty}
~\vspace{1.5cm}
\begin{center}{
\LARGE \bf
Quantum spin coverings and statistics
}
\end{center}
\vspace{1cm}
\centerline{\Large 
Ludwik D\c abrowski and Cesare Reina}
\vspace{5mm}
\begin{center}
{\it Scuola Internazionale Superiore Studi Avanzati,\\
Via Beirut 2-4, I-34014, Trieste, Italy}\\
\vspace{5mm}
\centerline{dabrow@sissa.it, reina@sissa.it}
\end{center}
\vspace{2cm}\noindent
\begin{abstract}\noindent
$SL_q(2)$ at primitive odd roots of unity $q^\ell =1$ is studied
as a quantum cover of the complex rotation group $SO(3,{\C})$,
in terms of the associated Hopf algebras of (quantum) polynomial functions. 
We work out the irreducible corepresentations,
the decomposition of their tensor products
and a coquasitriangular structure,
with the associated braiding (or statistics).
As an example, the case $\ell\! = \! 3$ is discussed in detail.
\end{abstract}
\vfill\eject
\section{Introduction}

As is well known, the twofold spin covering ~$\Z_2 \to Spin(2) \to SO(2)$~ 
is not universal and there are other (nontrivial) coverings with kernel $\Z_n$
or $\Z = \pi_1(SO(2))$ (for the universal one).
They are responsible for such features as the fractional or continuous spin 
and the associated anyonic statistics (see eg. \cite{W}, \cite{GH}).
In dimension three (or greater) $\pi_1(SO(3))=\Z_2$ 
and the spin covering ~$\Z_2 \to Spin(3) \to SO(3)$~ is universal, 
hence the only (projective) representations are of halfinteger spin 
(in addition to those of integer spin which are bona fide representations of $SO(3)$).
The quantum groups offer a possibility to refine this classification.
There is indeed a candidate for such a cover, 
the quantum group $SL_q(2)$ at the roots of unity \cite{C96}.
The case of the third root $q=e^{\frac{2\pi i}{3}}$ has been worked out
quite extensively mainly in terms of the quantum universal enveloping algebra
$U_q(sl(2))$ (cf. \cite{Ka}, \cite{CGT} and references therein) 
with the perspective to link to the Connes' algebra for the Standard Model. 
Also the Hopf algebra $A(SL_q(2))$ of `polynomials on $SL_q(2)$' 
at odd roots of unity $q^\ell =1$ has been studied.
We adopt this `quantum function' point of view since it is better suited 
to capture such topological properties as quantum (finite) covers,
which have to be introduced by hand when working in the universal 
enveloping algebra language. 
In Section 2, after recalling the essentials on $A(SL_q(2))$,
we shall present the coverings of $Spin(3,\C )$ and of $SO(3,\C )$.

The main point we are interested in this paper is to study further this 
new quantum symmetry and the braiding (statistics) of its corepresentations.
The labeling of irreducible corepresentations refines the notion spin. 
On the dual level, the irreducible representations of $U_q(sl(2))$ at roots of unity
are well known and the interesting decomposition of their tensor products
has been studied \cite{Ar}, \cite{A-GGS}, \cite{Ke}, \cite{PS}.
They are not of immediate use for us since the duality between  
$A(SL_q(2))$ and $U_q(sl(2))$ degenerates at roots of unity.
Actually it descends to a duality between some finite dimensional quotients
$A(F)$ and ${\bar {U}}_q$ of these Hopf algebras.
The irreducible representations of ${\bar {U}}_q$ 
are also known \cite{K} and by the (nondegenerate) duality 
they correspond to irreducible corepresentations of $A(F)$,
but this is not tantamount to those of $A(SL_q(2))$.
There exist another Hopf algebra of divided powers 
\cite{L1}, \cite{L2}, or a version of it due to  \cite{D-CK},
which is dual to $A(SL_q(2))$ \cite{D-CL}, \cite{D-C}.
Thus the irreducible corepresentations of $A(SL_q(2))$ should be 
in a correspondence with the irreducible representations 
of divided powers Hopf algebra, which are known.
This reasoning requires however some algebraic subtleties which we avoid.
For the sake of physicists' community, 
in Section 3 we provide a direct computational proof of the irreducibility.
The first nontrivial case $\ell\! = \! 3$ is discussed in Sect 3.1,
where we study the decomposition of tensor products 
of the irreducible corepresentations of $A(SL_q(2))$
and the question if it is possible to built a fundamental spinor
out of three `fractional' corepresentations.

As far as the braiding is concerned it is associated to 
a coquasitriangular structure.
In Section 4 we show that the coquasitriangular structure obtained from the 
standard  universal $R$ matrix associated to $U_q(sl(2))$ 
leads to a fairly exotic braiding, 
which however is consistent with the Bose-Fermi statistics 
of the usual (half)integer spin corepresentations.
As the first nontrivial example, the braiding in the case $\ell\! = \! 3$ 
is discussed in detail in Section 4.1.

Section 5 contains final remarks and conclusions.\\
\baselineskip15pt
\section{Preliminaries}

We start recalling the basic definitions and our notational conventions.

\subsection{The quantum group $SL_q(2)$}

Recall that for ${\mathbb C}\ni q \neq 0$, $A(SL_q(2))$ 
is the unital free algebra generated by $a,b,c,d$ over $\C$
modulo the ideal generated by the commutation relations
$$ab=qba\quad ac=qca,\quad bd=qdb,\quad bc=cb,\quad cd=qdc,\quad
 ad-da=(q-q^{-1})bc$$
and by the $q$-determinant relation $ad-qbc=1.$
The Hopf algebra structure is given by 
the following comultiplication $\Delta$, counit $\varepsilon$, and antipode $S$
defined on the generators (arranged as a $2\times 2$ matrix) by
\bea
\label{hs}
&
\Delta {\left( \ba{cc} a & b \\ c & d \ea \right)} =
{\left( \ba{cc} a & b \\ c & d \ea \right)} \ot 
{\left( \ba{cc} a & b \\ c & d \ea \right)}~,  \\
&
\varepsilon {\left( \ba{cc} a & b \\ c & d \ea \right)}=
{\left( \ba{cc} 1 & 0 \\ 0 & 1 \ea \right)}\, ,\;\;
S{\left( \ba{cc} a & b \\ c & d \ea \right)}=
{\left( \ba{cc} d &-q^{-1} b \\ -qc & a \ea \right)}\, ,\nonumber
\eea
where on the r.h.s. of the first equation the `line by columns' 
tensor product is understood. 
As a complex vector space, $A(SL_q(2))$ has a basis 
$a^ib^jc^k$, with $i,j,k\in\N$ ~and~ $b^ic^jd^k$ with $i,j\in\N, k\in\Z_+ $.
(We denote $\Z_+ =\{1, 2, 3\dots\}$ and $\N =\{0, 1, 2\dots\}$.)

\subsection{Finite quantum subgroups $F$ and ${\hat F}$}

From now on, unless stated differently, 
we set the parameter $q$ to be a (primitive) $\ell$-th root of unity
$\lambda =e^{\frac{2\pi i}{\ell}}$, for odd $\ell \geq 3 $.
We introduce two finite dimensional Hopf algebras of `functions
on finite quantum subgroups' $F$ and ${\hat F}$ of $SL_q(2)$,
needed in the sequel.

The Hopf algebra $A(F)$ is defined as the quotient Hopf algebra 
of $A(SL_q(2))$ modulo the ideal generated by the relations
\beq\label{F} 
a^\ell = 1 = d^\ell, ~~ b^\ell = 0 = c^\ell . 
\eeq
Let  $\pi_F$ denote the canonical projection, and $\tilde{t} := \pi_F(t)$.

We give now some information on $F$ 
(see \cite{DHS}, \cite{DNS} for the case $\ell =3$).
\begin{proposition}
\label{p1}
$A(F)$ satisfies the following properties:\\
i) as a complex vector space $A(F)$ is $\ell^3$-dimensional 
and its basis can be chosen as 
 $\tilde{a}^p \tilde{b}^r\tilde{c}^s$, where $p,r,s\in\{0,1,\dots, \ell\! -\! 1 \}$.\\
ii) $A(F)$ has a faithful representation $\varrho$
\beq \varrho(\tilde{a}) =  {\bf J} \ot {\bf 1}\sb\ell \ot {\bf 1}\sb\ell~, \quad
      \varrho(\tilde{b})  =  {\bf Q} \ot {\bf N} \ot {\bf 1}\sb\ell~, \quad
\varrho(\tilde{c})  = {\bf Q} \ot {\bf 1}\sb\ell \ot {\bf N}~, \eeq
 where for $i,j\in \{1,2,\dots, \ell\}$
\beq 
{\bf J}_{i,j} =  \left\{ \ba{cc} 
1 & \mbox{if $i\! =\! j\! +\! 1$ mod $\ell$}\\
0 & \mbox{otherwise} \ea , \right.\
     {\bf Q}_{i,j} = \left\{ \ba{cc} 
q^{-i} & \mbox{if $i\! =\! j$ }\\
0 & \mbox{otherwise} \ea ,\right.\
{\bf N}_{i,j} = \left\{ \ba{cc} 
1 & \mbox{if $i\! =\! j\! +\! 1$ }\\
0 & \mbox{otherwise} \ea . \right.
\label{rep}
\eeq 
iii) $A(F)$ has the `reduced' quantum plane as a quantum homogeneous space
(i.e. the algebra generated by $x$ and $y$ modulo the ideal generated by the relations 
$xy = q yx$, $x\sp\ell = 1$ and $y\sp\ell = 1$) which is isomorphic to 
$\M (\ell,\C)$.\\
iv) $F$ has a classical subgroup, defined as the group of characters of $A(F)$,
which is easily seen to be isomorphic to $\Z_\ell$.
Namely, for $i\in\{1, 2, \dots, \ell\}$ we have $\chi_i$
defined by their action on the generators as
$\chi_i (\tilde{a}) =q^i$, $\chi_i (\tilde{b}) =0$, $\chi_i (\tilde{c}) =0$ and 
$\chi_i (\tilde{d}) =q^{-i}$.
The Hopf algebra $A(\Z_\ell)$ appears as a quotient 
of $A(F)$ by the ideal generated by $\tilde{b}$, $\tilde{c}$
(which is also the intersection of the kernels of the characters).
\end{proposition}

Quite similarly, we define $A({\hat F})$
as the $2\ell^3$ dimensional quotient of $A(SL_q(2))$ 
modulo the relations 
\beq\label{barF} 
a^{2\ell} = 1 = d^{2\ell}, ~~ b^\ell = 0 = c^\ell . 
\eeq
Notice that the classical subgroup of ${\hat F}$ is 
the cyclic group $\Z_{2\ell}$. This group
`combines' the cyclic subgroup $\Z_{\ell}$ of $F$ 
with $\Z_{2}$ appearing in the classical spin cover.
In fact one has the exact sequence of groups
\be\label{zeta} 
0 {\lra} \Z_{2} {\lra} \Z_{2\ell} {\lra} \Z_{\ell} {\lra} 0 
\ee
which extends $\Z_{\ell} $ by the kernel $\Z_{2} $ of the classical
spin cover.
Notice that this extension is actually a semidirect (and not direct)
product iff $\ell$ is odd.

\subsection{Quantum group covering of $SL(2)$}

The Hopf subalgebra of $A(SL_q(2))$
generated by the $\ell^{th}$ powers
$$\alpha = a^l,\; \beta=b^l,\; \gamma= c^l,\; \delta=d^l$$
is isomorphic to the (commutative) Hopf algebra 
$$A(SL(2))={\mathbb C}[\alpha ,\beta ,\gamma , \delta ]\; /
<\alpha\delta - \beta\gamma -1>$$
with the restricted coproduct, counit and coinverse.
It is just the subalgebra of coinvariants of the coaction of $A(F)$ 
(as a quotient Hopf algebra).
It is known \cite{A} (see also \cite{DHS}, for the case $\ell =3$) 
that 
\begin{proposition}
\label{p2}
The sequence of algebras 
\be
A(SL(2)){\lra} A(SL_q(2)) \st{\pi\sb F}{\lra} A(F)
\label{cover}
\ee
is:\\
i) a (right, faithfully flat) Hopf-Galois extension of $A(SL(2)$ by $A(F)$
(quantum principal fibre bundle)\\ 
ii) a principal homogeneous Hopf-Galois extension 
(i.e., $A(F)$ is a quotient of the Hopf algebra $A(SL_q(2))$ 
by a Hopf ideal and $\pi\sb F$ is the canonical surjection)\\
iii) strictly exact (quantum quotient group)
\end{proposition}
(Notice that $\mbox{iii}\Rightarrow\mbox{ii}\Rightarrow\mbox{i}$, see
\cite{T}, \cite{PW}, \cite{S} for the relevant definitions).

Therefore (\ref{cover}) is a good candidate for a quantum covering 
of the spin group.
To be fully worth of this name, it should be checked to be nontrivial.
It can be seen that it is not totally trivial in the sense
that $A(SL_q(2))$ is not isomorphic to $A(SL(2))\ot A(F)$.
Another accepted notion to substitute the triviality 
for quantum principal bundles is however that of cleftness (or crossed products),
 c.f. e.g. \cite{DHS}.
It is not yet known if (\ref{cover}) is cleft and it seems to be a tough problem,
which is not tractable by the usual means (e.g. the theory of quantum
characteristic classes). A weaker result affirms that $A(SL_q(2))$ as a module
over $A(SL(2))$, which is finitely generated and projective (c.f. \cite{D-CL}),
is actually free \cite{DRZ}.
(The associated coherent sheaf of rank $l^3$ is free 
and the corresponding vector bundle $F$ over $SL(2)$ turns out to be trivial).
Moreover a set of $l^3$ generators can be chosen as \cite{DRZ}\\
 $$a^mb^nc^{s'},\qquad b^nc^{s''}d^r,$$
with the integers $m, n, r, s', s''$ in the range
$m\in\{1,\dots, \ell\! -\! 1\}$, $n, r\in\{0,\dots, \ell\! -\! 1\}$,
$s'\in\{m,\dots, \ell\! -\! 1\}$ and $s''\in\{0,\dots, \ell\! -\! r -\! 1\}$.

We expect that the quantum principal bundle (\ref{cover}) is actually 
non cleft and propose to employ it for a quantum spin covering 
in the next section.

\subsection{Quantum group covering of $SO(3,\C )$}

The (complex) group $SL(2)$ is isomorphic to the spin group $Spin(3,\C )$ 
and thus provides a twofold covering of the (complex) rotation group
$SO(3,\C )$.
This classical spin covering can be combined 
with the covering (\ref{cover}) as follows.
The Hopf algebra $A(SO(3,\C ))$ can be identified with the subalgebra
of even polynomials $A^+(SL(2))$ in $A(SL(2))$.
It can be seen that $A^+(SL(2))$ coincides with the 
coinvariants of the coaction of Hopf algebra $A({\hat F})$.
Let $\pi_{{\hat F}}$ denote the canonical projection and 
${\hat t} := \pi_{{\hat F}}(t)$.
In analogy with the proof of Prop. \ref{p2} it can be shown that
\begin{proposition}
\label{p3}
The sequence of Hopf algebras 
\be
A(SO(3,\C )) {\lra} A(SL_q(2)) \st{\pi\sb{{\hat F}}}{\lra} A({\hat F})\ ,
\label{lcover}
\ee
possesses the same nice properties i)-iii) as the sequence (\ref{cover}).
\end{proposition}
In particular (\ref{lcover}) is a quantum principal bundle and
referring to our remarks at the end of the previous subsection,
we mention that it is very likely non cleft.
In fact the relevant question about (\ref{lcover}) is if it is reducible 
to the subgroup $\Z_2$.
We expect that it is not the case, consistently with our conjecture 
about the noncleftness of (\ref{cover})
and thus propose 

\vspace{.2cm}\noindent
{\bf Definition}
~For any odd $\ell$, $\ell \geq 3 $, 
with $q =e^{\frac{2\pi i}{\ell}}$,
we call the principal fibre bundle (\ref{lcover})
quantum spin covering of the complex rotation group.

\section{Irreducible corepresentations}

There are two natural series of corepresentations of $A(SL_q(2))$.
The first one comes by restricting the coproduct to 
$W_n = \C\{\alpha^n, \alpha^{n-1}\gamma,\dots, \gamma^{n} \}$,
i.e. the span of monomials of degree $n$ in $\alpha = a^\ell$ and 
$\gamma = c^\ell$. The corepresentation $W_n$ is a `push forward' 
of the $\! n+\! 1$ dimensional (spin $n/2$) 
corepresentations of $A(SL(2))$ and thus obviously irreducible for all $n\in\N$.
The second one comes by restricting the coproduct to 
$Y_m ={\C}\{a^m, a^{m-1}c,\dots, c^m\}$,
i.e. the span of monomials of degree $m$ in $a, c$. 
More explicitly,
\bea\label{ym}
\Delta a^{m-h}c^h =&\sum_{r=1}^{m-h}\sum_{s=1}^h
q^{-r(h-s)}
{m-h\choose r}_{q^{-2}}{h\choose s}_{q^{-2}}
a^{m-h-r} b^r c^{h-s} d^s\otimes 
a^{m-(r+s)} c^{r+s} \nonumber\\
                    =&\sum_{k=0}^{m}\left(\sum_{r+s=k}
q^{-r(h-s)}
{m-h\choose r}_{q^{-2}} {h\choose s}_{q^{-2}}
a^{m-h-r} b^r c^{h-s} d^s\right) \otimes a^{m-k} c^k
\ ,
\eea
where 
$${k\choose j }_p = \frac{(k)!_p}{(k-j)!_p (j)!_p}, ~~~
(k)!_p = (k)_p (k-1)_p\dots(2)_p ~~{\rm and}~~ (k)_p=1+p+\dots +p^{k-1}\ .$$
It is indecomposable but not irreducible in general.
As we shall see, for $m\in\{0,1,\dots, \ell\! - \! 1\}$ $Y_m$ 
is indeed irreducible and we shall denote it $V_m$.
Also, for $m=n\ell\! - \! 1$, $n \in \Z_+$, 
it is irreducible but in fact equivalent to ${W_{n-1}}\ot V_{\ell-1}$.
More generally the corepresentations of the form
${W_{n}}\ot V_{m}$, with
$n\in {\mathbb N}, ~m\in \{0, 1, \dots , \ell -1\}$
are all irreducible as can be inferred
from \cite{D-CL}, \cite{D-C}, establishing the duality with a 
version \cite{D-CK} of divided powers algebra
\cite{L1}, \cite{L2} of which 
the classification of irreducible representations is known \cite{L1}.
Although straightforward, here we provide a direct computational proof. 
\begin{proposition}
\label{p4}
Set $m=m_0+\ell m_1$, with $0\leq m_0 \leq \ell-1,\; m_1\geq 0$.\\
a) For $m_1=0$ 
the comodule $V_{m_0} := Y_{m_0}$ is irreducible.\\
a$^\prime$ ) For every $m_1> 0$, the comodule $Y_{\ell-1+\ell m_1}$
(i.e. when $m_0=\ell -1$) is irreducible as well
and it is isomorphic to ${W_{m_1}}\ot V_{\ell-1}$.\\
b) when $0\leq m_0\leq \ell -2$ and $m_1\geq 1$, 
$Y_{m_0+\ell m_1}$ has a (maximal) subcomodule
isomorphic to $W_{m_1}\ot V_{m_0}$.
The quotient comodule is 
irreducible and isomorphic to $W_{m_1-1}\otimes V_{\ell -2-m_0}$.
\end{proposition}
{\bf Corollary.}~~
{\it 
The corepresentations
$W_n\ot V_m$ are irreducible for all
$n\in \N$ and $m\in\{0,1,\dots, \ell\! - \! 1\}$.
}

\vskip.2cm\noindent
{\it Proof of Proposition (\ref{p4}).~~}
The classical argument working for $q=1$ can be directly extended when $q$ 
is considered as an indeterminate and runs as follows.
Given a corepresentation $\rho$ of a Hopf algebra $A$ on a comodule $U$,
let $\rho_i^j$ be a matrix of elements of $A$ such that
$$u_i\mapsto \rho_i^j\otimes u_j ,$$
with respect to a basis $u_i$ ($i=1,\dots ,n)$ of  $U$.
There exist a coinvariant subcomodule $U'\subset U$ 
(with $\dim U'=k$, say) iff up to a conjugation by an invertible matrix $Z$ 
with elements in $\C [q,q^{-1}]$ the matrix $\rho$ takes a lower echelon form, 
i.e. iff 
$$\left(\matrix{\tau_1 &0\cr
                \tau_3 &\tau_4\cr}\right)
  \left(\matrix{z_1&z_2\cr z_3&z_4\cr}\right)=
  \left(\matrix{z_1&z_2\cr z_3&z_4\cr}\right)
  \left(\matrix{\rho_1&\rho_2\cr 
\rho_3&\rho_4\cr}\right)$$
where $U'$ is the span of the first $k$ elements of the transformed basis
and the block decomposition is given by the splitting $U=U'\oplus U/U'$.
In particular this requires that
$$\tau_1\left(\matrix{z_1&z_2\cr}\right)=
     \left(\matrix{z_1\rho_1\! +\! z_2\rho_3&
          z_1\rho_2\! +\! z_2\rho_4\cr}\right).$$
Notice that the $k\times n$ matrix $\left(\matrix{z_1&z_2\cr}\right)$ 
has rank $k$. Let $M$ be an invertible $k\times k$ submatrix.
We can write 
$$\tau_1 M=S ,\qquad \tau_1 M'=S'\ ,$$
where $S$ is the submatrix of the r.h.s. above corresponding to the 
columns of $M$ in $\left(\matrix{z_1&z_2\cr}\right)$
and prime denotes the submatrices with the complementary columns. 
Substituting, we get $k\times (n-k)$ linear relations over $\C [q,q^{-1}]$
among the elements $\rho_i^j$
$$S'=SM^{-1}M' \ .$$ 
Now let's have a closer look at the comodules $Y_m$.
The matrix elements of (\ref{ym}) are linear combinations of monomials
of degree $m$ in the generators $a ,b ,c ,d$.
For generic $q$, all $(m+3)(m+2)(m+1)/6$ of them appear in the 
sum on the r.h.s. of (\ref{ym}) and every monomial appears exactly in a single
matrix element i.e. two different matrix elements contain different monomials.
Therefore they are all  
linearly independent and, by the argument above, the $Y_m$ are irreducible. 

Recall that when $q =\lambda$ is a $\ell^{th}$ root of unity, $\lambda^\ell =1$,
the subalgebra generated by 
$\alpha =a^\ell,\beta =b^\ell,\gamma =c^\ell,\delta =d^\ell$
is central and isomorphic to the classical Hopf algebra $A(SL(2))$. 
Now several $q^{-2}$-binomial coefficients actually vanish 
when evaluated at $\lambda$. 
A simple way to control this is to use the fact that
the coproduct is an algebra homomorphism;
$$\Delta a^m=\Delta \alpha^{m_1}\Delta a^{m_0}\ ,$$
where $m=m_0+\ell m_1$. Hence
$$\sum_{r=0}^m {m\choose r}_{\lambda^{-2}} a^{m-r}b^r\otimes 
a^{m-r}c^r=$$
$$\sum_{i=0}^{m_1} \sum_{j=0}^{m_0}
{m_0\choose  j}_{\lambda^{-2}}{m_1\choose  i}_1
\alpha^{m_1-i}\beta^ia^{m_0-j}b^j\otimes\alpha^{m_1-i}\gamma^i 
a^{m_0-j}c^j,$$
giving the factorization formula (c.f. \cite{L1})
\be
\label{CHECK}
{m\choose r}_{\lambda^{-2}} = {m_0\choose r_0}_{\lambda^{-2}} {m_1\choose r_1},
\ee
where $m=m_0+\ell m_1,\; r=r_0+\ell r_1,\; 
0\leq m_0,r_0 \leq \ell-1,\; m_1,r_1\geq 0$
and last factor on the r.h.s. is just the ordinary binomial coefficient.
In particular all the binomial coefficients ${m\choose r}_{\lambda^{-2}}$ with 
$r_0>m_0$ vanish. 

When  $0\leq m = m_0\leq \ell -1$, this can not occur 
and the standard argument above yields the point a). 
For larger $m$'s however the comodule $Y_m$ is no more irreducible.  
Using again the  homomorphism property of the coproduct, we directly compute 
$$\Delta (a^{m-h}c^h) ~~=~~ \Delta (\alpha^{(m-h)_1}\gamma^{h_1})\times $$
$$\sum_{j=0}^{(m-h)_0}\sum_{t=0}^{h_0}\lambda^{-j(h_0-t)}
  {(m\! -\! h)_0\choose  j}_{\lambda^{-2}}{h_0\choose t}_{\lambda^{-2}}
  a^{(m-h)_0-j}b^j c^{h_0-t}d^t\otimes
a^{(m-h)_0+h_0-(j+t)}c^{j+t}\ .$$
\vspace{.2cm}\noindent
Notice that 
$$
m-h=
\left\{ \ba{cc} 
\ell (m_1-h_1)+(m_0-h_0) &\mbox{if ~$0\leq h_0\leq m_0$}\\
&\\
\ell(m_1-h_1-1)+(\ell+m_0-h_0) &\mbox{if ~$m_0+1\leq h_0\leq \ell -1$}\nonumber
 \ea , \right.\
\nonumber
$$
\vspace{.2cm}\noindent
and the formula above for $h_0\leq m_0$ reads 
\be
\label{h<m}
\Delta (a^{m-h}c^h) = 
\Delta (\alpha^{m_1-h_1}\gamma^{h_1}) \times
\ee
$$
\sum_{j=0}^{m_0-h_0}\sum_{t=0}^{h_0}\lambda^{-j(h_0-t)}
  {m_0\! -\! h_0\choose j}_{\lambda^{-2}}{h_0\choose t}_{\lambda^{-2}}
  a^{m_0-h_0-j}b^j c^{h_0-t}d^t\otimes
a^{m_0-(j+t)}c^{j+t}.
$$
This sum contains only monomials in $a,c$ of degree $m_0$.
So $W_{m_1}\ot V_{m_0}$ is an irreducible subcomodule.
If $m_0 =\ell -1$ this is isomorphic to the whole of 
$Y_{\ell -1 + \ell m_1}$.
This proves the point a$^\prime$) and the first statement of the point b).
To complete the proof of b) notice that for $m_0+1\leq h_0\leq \ell-1$ 
we have
$$
\Delta (a^{m-h}c^h) = \Delta (\alpha^{m_1-h_1-1}\gamma^{h_1})\times
$$
$$
  \sum_{j=0}^{\ell+m_0-h_0}\sum_{t=0}^{h_0}\lambda^{-j(h_0-t)}
  {\ell\! +\! m_0\! -\! h_0\choose  j}_{\lambda^{-2}} {h_0\choose t}_{\lambda^{-2}}
  a^{\ell+m_0-h_0-j}b^j c^{h_0-t}d^t\otimes a^{\ell+m_0-(j+t)}c^{j+t}.$$
Restricting the sum to  $j+t\leq m_0$ we can factor $a^\ell =\alpha$, while
restricting to $j+t\geq \ell$ we can factor $c^\ell=\gamma$, thus compensating
the $-1$ occurring in the exponent of the classical part of the coproduct 
and leaving in these two partial sums only monomials of degree $m_0$ in $a,c$. 
This cannot be done in the partial sum for $m_0+1\leq j+t\leq l-1$,
which gives
$$
\sum_{k=m_0}^{l-1}\sum_{s=0}^k \dots a^{l+m_0-h_0-k+s}b^{k-s} 
c^{h_0-s}d^s\otimes
a^{l+m_0-k}c^k=$$
$$\sum_{k'=0}^{l- m_0-2}\sum_s^{m_0 + 1 + k'} 
\dots a^{m'_0-h'_0-k'+s}b^{m_0+1+k'-s} 
c^{m_0+1+h'_0-t}d^t\otimes
(ac)^{m_0+1}a^{m'_0-k'}c^{k'}=$$
$$\lambda^n [(bc)^{m_0+1}\otimes (ac)^{m_0+1}] \Delta 
a^{m'_0-h'_0}c^{h'_0},$$
where $m'_0=l-m_0-2,\; h'_0=h_0+m_0+1$.
This expression contains monomials of degree larger than $m_0$.
Thus the quotient $Y_m / W_{m_1}\ot V_{m_0}$ is isomorphic to
$W_{m_1-1}\otimes V_{\ell -2-m_0}$ which completes the proof.

\vskip.2cm
We remark that the proof for generic $q$ is just the q-analogue 
of what happens in the classical case.
It is enough to notice that replacing the ordinary binomial coefficients
the q-analogue one gets (\ref{ym}) up to some nonvanishing factors.
This is why the corepresentation theory for generic $q$ is the same as 
the classical one.
In particular the argument above yields that the Clebsh-Gordan decomposition
holds as in the classical case
$$
Y_m\otimes Y_{m'} = Y_{m+m'} \oplus Y_{m+m'-2}\oplus \dots\ .
$$

In the case of $q$ being $\ell$-th root of unity it follows from Prop. \ref{p4}
that there is always at least one 
irreducible corepresentation of arbitrary dimension $D$,
and there are more (up to $\ell$) of them, according to how many
integers in $\{1,\dots, \ell\}$ divide $D$.

The decomposition rules 
of the tensor products $(W_n\ot V_m)\ot (W_{n'}\ot V_{m'})$ 
into irreducible corepresentations 
follow from the usual Clebsh-Gordan decomposition of $W_n\ot W_{n'}$ 
and the decomposition of $V_m\ot V_{m'}$, which obey a more complicated pattern.

\subsection{Decomposition of tensor products for $\ell =3$}

The previous discussion can be specified and simplified considerably 
in the simplest (nontrivial) case $\ell =3$. We have explicitly
the following matrices $\rho$ of three corepresentations $V_0$, $V_1$ and $V_2$
respectively
$$
1, ~
\left(\matrix{a&b\cr c&d\cr}\right)
~~{\mbox {and}} ~~
\left(\matrix{a^2&-q^2ab&b^2\cr 
ac&ad+q^{-1}bc&bd\cr c^2&-q^2cd&d^2}\right),
$$
as well as the usual form of $W_n$.

\noindent
According to the corollary of proposition \ref{p4} we see 
that there is one trivial (onedimensional, irreducible) corepresentation 
$V_0 = W_0$. 
In dimension two there are two 
(inequivalent) irreducible corepresentations $V_1$ and $W_1$. 
In dimension three there are also two: $V_2$ and $W_2$.
In dimension four $Y_3$ is indecomposable but not irreducible but there are 
two other irreducible corepresentations $W_3\ot V_0 = W_3$ and $W_1\ot V_1$.
In dimension five there is only one irreducible corepresentation $W_4$ 
($Y_4$ is indecomposable but not irreducible).
In dimension six there are three irreducible corepresentations:
$W_5$, $W_2\ot V_1$ and $W_1\ot V_2=V_5$.
A general pattern is that in any dimension $D$ there is always at least one 
irreducible corepresentation, if either $2$ or $3$ divide $d$ there are two
irreducible corepresentations and if $6$ divides $D$ there are three
irreducible corepresentations.

We give now the decomposition rules of the tensor products.\\	
Clearly, 
$$V_0\ot V_0 = V_0, ~~V_0\ot V_1 = V_1 ~{\rm and}~ V_0\ot V_2 = V_2\ .$$
Next, it can be seen that 
$$V_1\ot V_1 = V_0 \oplus V_2, ~V_1\ot V_2 = V_1 \oslash W_1 \oslash V_1
~{\rm and}~ V_2\ot V_2 = V_0 \oslash V_2 \oslash (W_1\oplus V_1)\oslash V_0\ ,$$
where $\oslash$ in an indecomposable corepresentation indicates that 
the left summand is a subcomodule while the right summand is a comodule 
only after quotienting the left one.
Noting that the tensor products in the opposite order decompose equivalently,
and recalling the usual decomposition
$W_n\ot W_{n'}= W_{|n-n'|} \oplus W_{|n-n'|+2}\oplus \dots \oplus W_{|n+n'|}$,
these rules permit to find a decomposition of tensor products of any number of 
$V_{m}\ot W_n$, with $n\in \N$ and $m\in\{1,2,3\}$.

A interesting question in this simplest (nontrivial) case of $\ell =3$
is if there is a possibility to built a fundamental fermion out of three anyons.
The decomposition rules permit to verify easily that the corepresentations 
$(V_1)^{\ot 3}$ and $(V_2)^{\ot 3}$ do not contain the fundamental spinor 
corepresentation $W_1$ as a subcorepresentation and the same is true 
for the tensor cube of the irreducible corepresentations $(V_m \ot W_n)$
if $m\in\{1,2\}$. They do however contain $W_1$ as a quotient (sub)corepresentation.
It is also worth to mention that the fundamental spinor $W_1$ subcorepresentation 
occurs nevertheless in e.g. the reducible but not decomposable representation $Y_3$ 
and thus also in its third tensor power $(Y_3)^{\ot 3}$.\\

\section{Braiding}

For general $q$ the quasitriangular structure on $U_q(sl(2))$ 
given by the well known universal element $R$ \cite{FRT} in 
(a suitable completion of) $U_q(sl(2))^{\ot 2}$, 
defines a coquasitriangular structure on $A(SL_q(2))$.
Its explicit form on the generators reads (c.f. \cite{K})
\be
{\cal R} \pmatrix{
a\ot a& a\ot b & a\ot c &a\ot d\cr
b\ot a& b\ot b & b\ot c &b\ot d\cr
c\ot a& c\ot b & c\ot c &c\ot d\cr
d\ot a& d\ot b & d\ot c &d\ot d
}
=
\pmatrix{
q^{-1/2} & 0 & 0        & q^{1/2}\cr
     0 & 0 & q^{-1/2}\! -\! q^{3/2} & 0\cr
     0 & 0 & 0        & 0\cr
     q^{1/2} & 0 & 0        &q^{-1/2}
}\ .
\ee
This structure provides a highly unusual (nonsymmetric and nondiagonal) braiding
of corepresentations $\rho$ and $\rho'$ of $A(SL_q(2))$
\be
\Psi (u_i \ot u_r') = \sum_{j,s} 
{\cal R} \left( {\rho'}_r^{~s} \otimes \rho_i^{~j}\right)
 u_s' \ot u_j \ .
\ee
In our situation, $q^\ell = 1$, it is not difficult however to verify that 
the corepresentations $W_n$ for $n$ odd 
(i.e. with halfinteger spin $n/2$)
are fermionic and the corepresentations $W_n$ for $n$ even  
(i.e. with integer spin $n/2$)
are bosonic, i.e. they obey 
\be
\Psi (w\ot w')= (-1)^{nn'} w'\ot w {\rm ~for~} w\in W_n, ~w'\in W_{n'}\ .
\ee
Thus the exotic coquasitriangular structure ${\cal R}$ 
passes an important consistency test of the agreement 
with the usual spin-statistics  relation in dimensions $d\geq 3$.
The braiding of $V$'s and $W$'s is also quite simple:
\be
\Psi (v\ot w)= (-1)^{mn} w\ot v   {\rm ~for~} v\in V_m, ~w\in W_{n} \ .
\ee
Instead the braiding of $V$'s among themselves is highly exotic
(even comparing with the anyonic one).
We report it in the Subsection 4.1 for the case $\ell =3$.
Clearly the tensor products $V_m\ot W_n$ carry the combined statistics
according to the usual hexagon conditions for $\Psi$ 
(see eg. \cite{M95}).

\subsection{Braiding in the case $\ell =3$}

As for general $\ell$, the braiding of the corepresentations 
$W_n$ with $W_n'$ is exactly the classical one, i.e. the trivial twist 
except when $nn'$ is odd when it is (-) the twist.
Also, $W_n$ have the trivial braiding with $V_0$ and with $V_2$ 
and (-) the twist with $V_1$.
The braiding of $V$'s among themselves is as follows:\\
The braiding of $V_1$ and $V_1$ reads:
\be\label{11}
{\Psi} \pmatrix{
a\ot a\cr
a\ot c\cr
c\ot a\cr
c\ot c
}
=
\pmatrix{
q^{-1/2} & 0       & 0                & 0\cr
       0 & 0       & q^{1/2}          & 0\cr
       0 & q^{1/2} & 1\! +\! q^{-1/2} & 0\cr
       0 & 0       & 0                & q^{-1/2}
}
\pmatrix{
a\ot a\cr
a\ot c\cr
c\ot a\cr
c\ot c
}
\ .
\ee
(Notice that $\Psi$ has a nonsiple tensor $a\ot c-q c\ot a$
as an eigenvector with eigenvalue 1
and the complex span of $ a\ot a, qa\ot c + c\ot a, c\ot c$ 
as an eigenspace with eigenvalue $q^{-1/2}$).\\
Next, the braiding of $V_1$ and $V_2$ reads:
\be\label{12}
{\Psi} 
\pmatrix{
a^2\ot a\cr
a^2\ot c\cr
ac\ot a\cr
ac\ot c \cr
c^2\ot a\cr
c^2\ot c
}
=
\pmatrix{
q^{2}  & 0 & 0 & 0        & 0   & 0\cr
     0 & 0 & 0 & q        & 0   & 0\cr
     0 & 1 & 0 & q^{2}\! -\! q & 0   & 0\cr
     0 & 0 & 0 & 0        & 1   & 0\cr
     0 & 0 & q & 0        & 1\! -\! q & 0\cr
     0 & 0 & 0 & 0        & 0   & q^{2}
}
\pmatrix{
a\ot a^2\cr
a\ot ac\cr
a\ot c^2\cr
c\ot a^2\cr
c\ot ac \cr
c\ot c^2\
}
\ .
\ee
The opposite braiding of $V_2$ and $V_1$ reads:
\be\label{21}
{\Psi} 
\pmatrix{
a\ot a^2\cr
a\ot ac\cr
a\ot c^2\cr
c\ot a^2\cr
c\ot ac \cr
c\ot c^2\
}
=
\pmatrix{
q^{2} & 0 & 0 & 0        & 0   & 0\cr
    0 & 0 & 1 & 0        & 0   & 0\cr
    0 & 0 & 0 & 0        & q   & 0\cr
    0 & q & q\! -\! q^2  & 0   & 0  & 0\cr
    0 & 0 & 0 & 1        & 1\! -\! q^2 & 0\cr
    0 & 0 & 0 & 0        & 0   & q^2
}
\pmatrix{
a^2\ot a\cr
a^2\ot c\cr
ac\ot a\cr
ac\ot c \cr
c^2\ot a\cr
c^2\ot c
}
\ .
\ee
Finally, the braiding of $V_2$ and $V_2$ reads:
\be\label{22}
{\Psi} 
\pmatrix{
a^2\ot a^2\cr
a^2\ot ac\cr
a^2\ot c^2\cr
ac\ot a^2\cr
ac\ot ac \cr
ac\ot c^2\cr
c^2\ot a^2\cr
c^2\ot ac\cr
c^2\ot c^2
}
=
\pmatrix{
q & 0 & 0   & 0         & 0 & 0 & 0              & 0         & 0\cr
0 & 0 & 0   & 1         & 0 & 0 & 0              & 0         & 0\cr
0 & 0 & 0   & 0         & 0 & 0 & q              & 0         & 0\cr
0 & 1 & 0   & 1\! -\! q & 0 & 0 & 0              & 0         & 0\cr
0 & 0 & 0   & 0         & 1 & 0 & 1\! -\! q^2    & 0         & 0\cr
0 & 0 & 0   & 0         & 0 & 0 & 0              & 1         & 0\cr
0 & 0 & q^2 & 0         & q\! -\! 1 & 0 & -(q\! -\! 1)^2 & 0         & 0\cr
0 & 0 & 0   & 0         & 0 & 1 & 0              & 1\! -\! q & 0\cr
0 & 0 & 0   & 0         & 0 & 0 & 0              & 0         & q^2
}
\pmatrix{
a^2\ot a^2\cr
a^2\ot ac\cr
a^2\ot c^2\cr
ac\ot a^2\cr
ac\ot ac \cr
ac\ot c^2\cr
c^2\ot a^2\cr
c^2\ot ac\cr
c^2\ot c^2
}
\ .
\ee
The resulting braiding of the ireeducible corepresentations 
$W_n\ot V_m$ and $W_{n'}\ot V_{m'}$ can be obtained from the above braidings
of $V$'s and $W$'s using the hexagon conditions for braiding.

\section{Final remarks}

We remark in a connection with the point iii) of Prop. \ref{p1} 
that $\M(3,\C)$ occurs as a direct summand of Connes' interior algebra ${\cal A}$
for the Standard Model \cite{C94}.
Also, the algebra $M(3,\C) \oplus M(2,\C) \oplus \C$, close to ${\cal A}$,
coincides with the semisimple part of the algebra $U_q(sl(2))$ 
at cubic roots of unity, which was extensively studied
(cf. \cite{CGT}, \cite{Ka}, \cite{DNS} and references therein).

We mention also some related works. In \cite{M} the braided group $B(SL_q(2))$
with the braiding induced by the universal $R$-matrix via the 
right adjoint corepresentation has been described.
In \cite{DMAP} the fractional supersymmetry has been discussed.
The interesting papers \cite{O1} and \cite{O2} investigate 
the spin-statistics relation in the supergroup framework and its 
bosonisation to an (ordinary) Hopf algebra.
In \cite{GM} the noncommutative cohomology and electromagnetism on 
$SL_q(2)$ at roots of unity have been studied.
 
We summarize our work by saying that in noncommutative
geometry the spin and the statistics 
in dimensions $\geq 3$ looks more similar to the case of dimension $=2$.
In fact, to the best of our knowledge, besides the parafermion case, 
it opens a possibility of unusual statistics in quantum theories 
in dimension $\geq 3$.
It is however quite encouraging that the corepresentations $W_n$
of the usual $Spin(3,\C )$ maintain their Bose-Fermi statistics
in agreement with the usual spin-statistics theorem 
in (local) relativistic quantum field theory.

There are some open problems regarding the noncommutative spin covers.
From the mathematical point of view certainly the 
question about the cleftness of (\ref{cover}) and the (non)reducibility
of (\ref{lcover}) to $\Z_2$ should be settled. 
Also the issue of reality or *-structure and the reductions 
to $SU(2)$ and $SO(3)$ are very important. 
Of course the covers of relativistic symmetries 
(Lorentz and Poincar\'e) should be then worked out as well.
Another task is to employ the quantum covers as structure groups for
bundles on (classical or quantum) spaces. 
These and related topics are currently under investigation.

From the physical point of view, besides the hypothetical relation
with quantum symmetries behind the Standard Model,
the main question concerns a possible role of the quantum covers of spin,
e.g. in local quantum field theory.
Indeed from the discussion above we see that the extension (\ref{zeta})
is actually nontrivial and hence there is an interesting mixing up
between the exotic statistics of $V$'s and spin.
This may be a key test for physical applications.
Also the use in physics of the indecomposable corepresentations,
their irreducible subcomodules and of their quotients require a further study.

\end{document}